\documentclass[11pt,leqno]{amsart}
\usepackage{amsmath}
\usepackage{amsthm}
\usepackage{amssymb}
\usepackage{latexsym}
\usepackage[dvips]{graphicx}
\theoremstyle{definition}

\theoremstyle{plain}

\begin{document}
\setcounter{page}{1}
\begin{flushleft}
\scriptsize{{Unpublished Research Note\/} {\bf 1}\,(October 2009),
\\
\copyright K. O. Babalola\\
Department of Mathematics\\
University of Ilorin\\
Ilorin, Nigeria\\
kobabalola\symbol{64}gmail.com\\
khayrah.babalola\symbol{64}gmail.com\\
abuqudduus\symbol{64}yahoo.com\\
babalola.ko\symbol{64}unilorin.edu.ng\\
http://unilorin.academia.edu/KunleOladejiBabalola/Teaching}
\end{flushleft}
\vspace{16mm}

\begin{center}
{\normalsize\bf An invitation to the theory of geometric functions} \\[12mm]
     {\sc K. O. BABALOLA} \\ [8mm]

\begin{minipage}{123mm}
{\small {\sc Abstract.}
This note is an invitation to the theory of geometric functions. The foundation techniques and some of the developments in the field are explained with the mindset that the audience is principally young researchers wishing to understand some basics. It begins with the basic terminologies and concepts, then a mention of some subjects of inquiry in univalent functions theory. Some of the most basic subfamilies of the family of univalent functions are mentioned. Main emphasy is on the important class of Caratheodory functions and their relations with the various classes of functions, especially the techniques for establishing results in those other classes when compared with the underlying Caratheodory functions. This is contained in Section 4. Examples based on this technique are given in the last section. Since the target audience is the uninitiated, the difficult proofs are not presented. The elementary proofs are explained in the simplest terms. Footnotes are made to further explain some not-immediately obvious points. The references are mostly standard texts. The interested may consult experts for the most recent references in addition to those contained in the cited texts. Hopefully, this may as well profit even the initiated who intends to research in this field.}
\end{minipage}
\end{center}

 \renewcommand{\thefootnote}{}

\def\iff{if and only if }
\def\S{Smarandache }
\newcommand{\norm}[1]{\left\Vert#1\right\Vert}

\vskip 12mm

\baselineskip=0.25in

\section{Introduction}
Let us begin by saying that: functions being studied in this subject area are generally complex-valued and analytic in a chosen domain. They may be of several variables. However, our focus in this note is largely on those functions which are of one complex variable. Such a function (say $g$) is said to be analytic (regular or holomorphic) at a point $z_0$ in its domain if its {\em derivative} exists there. Because these functions are analytic (and thus are continuously differentiable), they have Taylor series developments in their domain. They are thus expressible in certain series form with centres at (say) $z_0$. Since by simple translation the nonzero centres $z_0$ may be {\em shifted} to zero, we may assume without loss of generality that the centres of the series developments of these functions are the origin. Thus an analytic function $g$ may be expressed as:
$$g(z)=b_0+b_1z+b_2z^2+b_3z^3+\cdots.$$
The coefficient $b_k=g^{(k)}(0)/k!$ and is easily obtained from the Cauchy integral formula
$$g^{(k)}(z)=\frac{k!}{2\pi i}\int_\Gamma\frac{g(w)}{w-z}dw$$
where $\Gamma$ is a rectifiable simple closed curve containing $z$ and $g$ is analytic inside and on it.\vskip 2mm

{\bf The unit disk:} We would assume the domain of $g$ to be the unit disk $E=\{z:|z|<1\}$. Any justification for this? Yes. The Riemann mapping theorem guaranttees that any such domain ({\em simply connected}) in the complex plane can be mapped {\em conformally} to any other with similar description. Put differently, Riemann showed that there always exists an analytic function that maps one simply connected domain to another also with similar description. This epoch assertion of Riemann seemed to have lacked full flavour or strength, until the birth of the theory of univalent functions. In 1907, Koebe discovered that analytic and univalent mappings have the nice quality of the Riemann assertion \cite{OP}:\vskip 2mm

\begin{quote}
{\em If $z_0\in D$, then there exists a unique function $g$, analytic and univalent which maps $D$ onto the the open unit disk $E$ in such a way that $g(z_0)=0$ and $g'(z_0)>0$.} 
\end{quote}
\vskip 2mm

Thus with the univalence (and thus conformality) of $g$, heads need not ache regarding specifics of the geometry of any simply connected domain in the complex plane as many varieties of problems about such domains are invariably reducible to the special case of the open unit disk.\vskip 2mm  

{\bf Normalization:} The function $g$ is normalized such that:\vskip 2mm

\begin{itemize}
\item[(i)] it takes the value zero at the origin (that is it takes the origin to the origin, $g(0)=0$) and\vskip 2mm

\item[(ii)] its derivative takes the value 1 at the origin, that is $g'(0)=1$.\vskip 2mm
\end{itemize}

Why? Observe this from the Riemann assertion that, without loss of generality, we may take $z_0=0$ so that the assertion becomes:

\begin{quote}
{\em If $D$ contains the origin, then there exists a unique function $g$, analytic and univalent which maps $D$ onto the the open unit disk $E$ in such a way that $g(0)=0$ and $g'(0)>0$.} 
\end{quote}
\vskip 2mm

The requirement $g(0)=0$ and $g'(0)>0$ is exactly the reason for normalization. Now how is this to be achieved? Define 
$$f(z)=\frac{g(z)-b_0}{b_1}$$
provided the coefficient $b_1\neq 0$. Is this condition true of all analytic function $g$? Definitely not! The analytic function $g(z)=z^2$ is a counterexample. However, there are yet many others so normalizable. So, we know, sure, that the class of normalizable analytic functions is nonempty. Fortunately, there exists a subset of them which have a nice underlying property. Alas, these are those that are injective or one-to-one. In geometric functions' parlance, such functions are variously called {\em univalent, simple}, {\em schliht} (German) or {\em odnolistni} (Russian). They are functions which do not take on the same value twice. That is if $z_1, z_2$ are points in the domain (say $D$) of $g$, then 
$$g(z_1)=g(z_2)\Longrightarrow z_1=z_2.$$
Put in another way,
$$z_1\neq z_2\Longrightarrow g(z_1)\neq g(z_2),\;\;\text{for all}\;\;z_1,z_2\in D.$$
It is not so difficult to see graphically that $f$ is injective if and only if $f'(z)\neq 0^{(1)}$\footnote{ $^{(1)}$A function $f$ is one-to-one if and only if it does not turn in its domain, for it does, then in some neighborhood of its turning point, it must assign the same value twice. The mathematical presentation of the {\em never turning} property is: $f'(z)\neq 0$ for all $z\in D$. This is easily seen graphically on $\mathbb{R}$.}, that is it does not have zero $b_1$. In other words $f$ is injective if and only if it never turns in its domain. A simple analytic proof is that if by contradiction it is assumed that they do, then for sufficiently small $z$, $g$ may be approximated (taking o($z^3$) as zero) by:
$$g(z)\approx b_0+b_2z^2$$
in which case $g$ looses univalence.\vskip 2mm

Now we are guarantteed that with the univalence of $g$, the desired normalization can be effected and we thus isolate them and denote them by $S$, say. Furthermore we represent them by:
$$f(z)=z+a_2z^2+\cdots\eqno(1)$$
where $a_k=b_k/b_1$, $k=2,3,\cdots$ and $b_1\neq 0$.\vskip 2mm

{\bf The range of} $f$: Is the nomenclature {\em geometric function theory} a misnomer or not for this field of study? No, it isn't. In the words of Macgregor:\vskip 2mm

\begin{quote}
{\em The significance of geometric ideas and problems in complex analysis is what is suggested by the term geometric function theory. These ideas also occur in real analysis, but geometry has had a much greater impact in complex analysis and it is a very fundamental aspect of its vitality.}
\end{quote}
\vskip 2mm

Duren \cite{PL} adds:\vskip 2mm

\begin{quote}
{\em The interplay of geometry and analysis is perhaps the most fascinating aspect of complex function theory. The theory of univalent functions is concerned primarily with such relations between analytic structure and geometric behaviour.}
\end{quote}
\vskip 2mm

The ranges of these functions describe various nice geometries and classical characterizations. An example is: if $f$ is a normalized analytic and univalent function in $E$, then its range contains some disk $|w|<\delta$. Furthermore, the ranges of some of them describe {\em star, close-to-star, convex, close-to-convex or linearly accessible, spiral} geometries: some in certain directions, some uniformly, some with respect to conjugate symmetric points and so on. These functions whose ranges describe certain geometries are thus known as geometric functions. Furthermore, their study is also known as Geometric Functions Theory.\vskip 2mm

In particular, a region of the complex plane is said to have star geometry with respect to a fixed point in it if every other point of it is visible from the fixed point. In other words, a ray or line segment issuing from the fixed point inside it to any other point of it lies entirely in it. If a region has star geometry with respect to every point in it, it is called convex. That is, the line segment joining any two points of this region lies entirely inside it.\vskip 2mm   

Functions whose ranges have star geometry are known as star functions while those whose ranges have convex geometry are called convex. This same notion is expressed in many other classes of functions.\vskip 2mm
 
{\bf Between analysis and geometry}: Any connections? Yes. Researchers have made groundbreaking discoveries between analysis and geometry. They have succeeded not only in describing those geometries in succint mathematical terms, but also in establishing close links between certain prescribed properties of analytic functions and the geometries of their ranges.  For example, if a function $f$ maps the unit disk onto a star domain, then the real part of the quantity $zf'/f$ is positive. The converse is also true. Similarly, if $f$ maps the unit disk onto a convex domain, then the real part of the quantity $zf''/f'$ is greater than -1. The converse is also true. Furthermore, with the truth also of the converse, if $f$ maps the unit disk respectively onto a close-to-star, close-to-convex or linearly accessible, spiral domain, then the real part of the quantities $f/g$, $f'/g'$ and $e^{i\theta}zf'/f$ is positive, where $g$ is convex. These, perhaps, have led to the thinking that if the positivity condition of the real parts of many of these quantities is necessary as well as sufficient for univalence, then what can be said of other quantities such as $f'$, those involving higher derivatives or defined by certain operators and more recently of linear combinations of two or more of such quantities?
\vskip 2mm

{\bf Any examples?}: Yes. The leading member of the large family of univalent functions is the famous Koebe functions given by
$$k(z)=\frac{z}{(1-z)^2}=z+2z^2+3z^3+\cdots.$$
The Koebe function maps the open unit disk onto the entire complex plane except a slit along the negative real axis from $-\frac{1}{4}$ through to $-\infty$. For many problems regarding the entire family of univalent functions (and some subsets of it), the Koebe function assumes the best possible extremum. We demonstrate this by examples in latter sections. A trivial member of the family is the identity mapping $f(z)=z$. The identity mapping is ubiquitous; it can be found in any subclass of the class of univalent functions.\vskip 2mm

{\bf Best possible property}: The class of univalent functions and many subclasses of it are being studied in the abstract sense. Many characterizations of them apply in the general sense to all members of the class under consideration as is the case with many subjects of pure mathematics. Now, if a property or characterization $T$ on a class of functions (or any set, $J$, for that matter) is such that there exists a member of the class $J$ assuming the extremum, then such a property is said to be {\em best possible} on $J$. For example, in $\mathcal{S}$, the coefficient characterization inequality $|a_2|\leq 2$ is best possible since the Koebe function, $k(z)$, which is a member of $\mathcal{S}$, takes the equality. This is to say the property $|a_2|\leq 2$ cannot be made better as along as the Koebe function is a member of the set under consideration. The synonyms of ``best possible'' as can be found in usage by many workers in this field are ``{\em sharp}'' and ``{\em cannot be improved}''.

\section{Some subjects of inquiry}
A wide range of problems of mathematical analysis are being solved in the theory of geometric functions as many as their results are being applied in many branches of mathematics, physical sciences and engineering. Before long let us refer to the great compilation  by S. D. Bernardi:\vskip 2mm

\begin{quote}
{\em Bibliography of Schlicht functions}, Courant Institute of Mathematical Sciences, New York University, 1966; Part II, {\em ibid}, 1977. Reprinted with Part III added by Mariner Publishing Co. Tampa, Florida, 1983.,
\end{quote}
which itemize the many subject areas of the Geometric Functions Theory plus the list of the many research outputs in those areas.\vskip 2mm

We now begin our few mention of them by first noting the fact that these univalent functions exist infinitely in nature so much so that the simple definition, $f(z_1)=f(z_2)\Longrightarrow z_1=z_2$ or its equivalent $z_1\neq z_2\Longrightarrow f(z_1)\neq f(z_2)$, cannot be used in general to identify, isolate or recognize many of them. This has given birth to several new methods of mathematical analysis with this sole aim. In particular these methods came under what is usually refered to as:\vskip 2mm

{\bf Sufficient conditions for univalence}. Results in this direction are as many as there are researchers in the field. They continue to appear in prints with no end in sight. Notable and simplest among them is the statement:\vskip 2mm

\begin{quote}[Noshiro-Warschawski Theorem \cite{PD}]
{\em If $f$ is analytic in a domain $D$ and Re $f'(z)>0$ there, then $f$ is univalent there.}
\end{quote}
\vskip 2mm

The proof of the above univalence condition depends on the fact that the function $f$ is defined on a line segment joining any two distinct points of its domain, say, $L:tz_2+(1-t)z_1$, so that by the transformation $z=tz_2+(1-t)z_1$ ($dz=(z_2-z_1)dt)$, we have$^{(2)}$\footnote{ $^{(2)}$The linear segment $z:=tz_2+(1-t)z_1$ implies that when $z=z_2$ then $(1-t)z_2=(1-t)z_1$, which holds if and only if $t=1$ since $z_1\neq z_2$. Similarly when $z=z_1$ we have $tz_2=tz_1$, which holds also if and only if $t=0$ since $z_1\neq z_2$, thus leading to the new integral in the proof.}
$$\aligned
f(z_2)-f(z_1)
&=\int_{z_1}^{z_2}f'(z)dz\\
&=(z_2-z_1)\int_0^1f'(tz_2+(1-t)z_1)dt\neq 0\;\;\text{since Re} f'(z)>0. 
\endaligned$$

In fact, the assertion of the Noshiro-Warschawski theorem is contained in an equivalent but more general statement, which is:\vskip 2mm

\begin{quote}[Close-to-convexity \cite{PD}]
{\em If $f$ is analytic in a domain $D$ and if for some convex function $g$,  Re $f'(z)/g'(z)>0$ there, then $f$ is univalent there.}
\end{quote}
\vskip 2mm

\begin{proof}
Let $D$ be the range of $g$ and consider $h(w)=f(z)=f(g^{-1}(w))$, $w\in D$. Then
$$h'(w)=\frac{f'(g^{-1}(w))}{g'(g^{-1}(w))}=\frac{f'(z)}{g'(z)}$$
so that Re $h'(w)>0$ in $D$. Thus $h(w)=f(z)$ is univalent$^{(3)}$\footnote{ $^{(3)}$Since $h(w)=f(z)$, then $h'(w)dw=f'(z)dz$. But $z=g^{-1}(w)$, that is $w=g(z)$ so that $dw=g'(z)dz$. Hence $h'(w)dw=f'(z)dw/g'(z)$, that is $h'(w)=f'(z)/g'(z)$ as in the proof.}.
\end{proof}
\vskip 2mm

Perhaps more than any other, this subject has led to identifying many more subfamilies of the class of univalent functions in the unit disk. Some of these subclasses are discussed in Section 3.\vskip 2mm

Close to this is the inquisition about which transformations preserve univalence in the unit disk. The most basic ones being: conjugation, $\overline{f(\bar{z})}$; rotation, $e^{-i\theta}f(e^{i\theta}z)$; dilation, $f(rz)/r$ for $0<r<1$; disk automorphism, $[f((z+\sigma)/(1+\bar{\sigma}z))-f(\sigma)]/[(1-|\sigma|^2)f'(\sigma)]$, $\sigma\in E$; omitted-value, $\xi f(z)/[\xi-f(z)]$, $f(z)\neq\xi$, $\xi\in E$; square root, $\sqrt{f(z^2)}$; and the composition/range transformations, $\varphi(f(z))$ where $\varphi$ is similarly normalized analytic and univalent but in the range of $f$. All the transformations are easily verified via the definition $f(z_1)=f(z_2)\Longrightarrow z_1=z_2$, except the square root transformation, which requires a little explanation$^{(4)}$\footnote{ $^{(4)}$Note that the function $g(z)=\sqrt{f(z^2)}=z+c_3z^3+c_5z^5+\cdots$ is an odd analytic function such that $g(-z)=-g(z)$. So if $g(z_1)=g(z_2)$, then $f(z_1^2)=f(z_2^2)$ and thus $z_1^2=z_2^2$. That is $z_1=\pm z_2$. But if $z_1=-z_2$, then $g(z_1)=g(z_2)=g(-z_1)=-g(z_1)$, so that $2g(z_1)=0$ and $z_1=0$ since $f(0)=0$ only at the origin. Thus we have $g(z_1)=g(z_2)\Longrightarrow z_1=z_2$, which shows that $g$ is univalent.}. Advances in the subject have led to consideration of more difficult transformations, particularly those ones which are solutions of certain linear/nonlinear differential equations. The simplest form of this is what came to be know as the Libera integral transform defined as:
$$\mathcal{J}(f)=\frac{2}{z}\int_0^zf(t)dt,\;\text{(See \cite{RJ})}.\eqno{(2)}$$

The Libera integral is the solution of the first-order linear differetial equation: $zf'(z)+f(z)=2g(z)$. Various other integrals have been considered, many being generalizations of the Libera integral. Transformations of this type examine the nature and propeties of the solutions of certain differential equations given that $f$ has some known properties and or the extent of such properties being transferable to the solutions.\vskip 2mm

{\bf Radius problems}. If we suppose that some transformations or geometric conditions fail to preserve univalence (for instance) in the unit disk, then it is natural to ask if such transformations (or conditions) could preserve it in any subdisk $E_0=\{z:|z|<\rho<1\}\subset E$. Problems of this sort became known as radius problems. More precisely, it is about finding the radius $\rho$ of the largest subdisk $E_0$ in which certain transformations of a univalent function $f$ or some geometric conditions guarantees univalence. This radius $\rho$ is particularly known as the radius of univalence (for instance). By ``for instance'' we imply that this notion is not restricted to the subject of univalence only. In fact, and interestingly, this has raised many more questions like: the radius of starlikeness, convexity, close-to-convexity and many more. A basic result in this direction is:\vskip 2mm

\begin{quote}[Noshiro, Yamaguchi \cite{KY}]
{\em If $f$ satisfy Re $f(z)/z>0$ in $E$, then it is univalent in the subdisk $|z|<\sqrt{2}-1$}
\end{quote}
\vskip 2mm

{\bf Convolution or Hadamard product}. Let $f(z)=a_0+a_1z+a_2z^2+\cdots$ and $g(z)=b_0+b_1z+b_2z^2+\cdots$ be analytic funtions in the unit disk. The convolution (or Hadamard
product) of $f(z)$ and $g(z)$ (written as $(f*g)(z)$) is defined as
$$(f*g)(z)=z+\sum_{k=2}^\infty a_kb_kz^k.$$
The concept of convolution arose from the integral
$$h(r^2e^{i\theta})=(f*g)(r^2e^{i\theta})=\frac{1}{2\pi}\int_0^{2\pi}f(re^{i(\theta-t)})g(re^{it})dt,\;\;\;r<1$$
and has proved very resourceful in dealing with certain problems of the theory of analytic and univalent functions, especially closure of families of functions under certain transformations. This is since many a transformation of $f$ is expressible as convolution of $f$ with some other analytic function, sometimes with predetermined behaviour. It is natural, therefore, to desire to investigate the convolution properties of many classes of functions. For example the Libera transform (2) is the convolution $\mathcal{J}=g\ast f$ where $g$ is the analytic function
$$g(z)=z+\sum_{k=2}^\infty\frac{2}{k+1}z^k.$$ 
This function $g$ has some nice geometric properties which may pass on to the Libera transform via the convolution as would be found in literatures through further studies.\vskip 2mm

{\bf Coefficient inequalities}. A close look at the series development of $f$ suggests that many properties of it like the growth, distortion and in fact univalence, may be affected (or be told) by the size of its coefficients. Duren says:\vskip 2mm

\begin{quote}
{\em In most general form, the {\em coefficient problem} is to determine the region of $\mathbb{C}^{n-1}$ occupied by the points $(a_2,\cdots,a_n)$ for all $f\in S$. The deduction of such precise analytic information from the geometric hypothesis of univalence is exceedingly difficult.}
\end{quote}
\vskip 2mm

The most contained in this part of this article are sourced from the survey by Duren \cite{PL}, which is ample for detailed issues regarding the coefficient problems in the field.\vskip 2mm

The coefficient problem has been reformulated in the more special manner of estimating $|a_n|$, the modulus of the $n$th coefficient. Perhaps, no problem of the field has challenged its people as much as the coefficient problem. As early as in 1916, Bieberbach conjectured that the $n$th coefficient of a univalent function is less or equal to that of the Koebe function. In mathematical language, he says:\vskip 2mm

\begin{quote}
{\em For each function $f\in S$, $|a_n|\leq n$ for $n=2,3,\cdots$. Strict inequality holds for all $n$ unless $f$ is the Koebe function or one of its rotations.}
\end{quote}

The conjecturer, Bieberbach, himself proved that $|a_2|\leq 2$ as a simple corollary to the {\em area theorem}$^{(5)}$\footnote{ $^{(5)}$See \cite{PD}, page 29 for the {\em area theorem} and the proof of Bieberbach theorem ($|a_2|\leq 2$).}, which is due to Gromwall. The third was settled in 1923 by Loewner. The fourth was solved in 1955 by Garabedian and Schiffer, while in 1960 Charzynski and Schiffer gave an elementary proof of same result. The proofs for the fifth and sixth came several years latter. Thereafter, the great puzle had remained unsolved until only recently when, precisely 1985, De Brange announced the final solution to the notorious conjecture. In total, the conjecture had stood for sixty-nine years unsolved! These long years were not unproductive however, as the conjecture had inspired the development of important new methods and techniques in the theory in particular and complex analysis in general.\vskip 2mm

Closely related to the Bieberbach conjecture is that of finding the sharp estimate for the coefficients of odd univalent functions, which has the most general form of the square root transformation of a function $f\in\mathcal{S}$:
$$l(z)=\sqrt{f(z^2)}=z+c_3z^3+c_5z^5+\cdots.$$
For odd univalent functions, Littlewood and Parley in 1932 proved that for each $n$ the modulus $|c_n|$ is less than an absolute constant $A$, (which their method showed is less than 14) and they added the footnote ``No doubt the true bound is given by $A=1$'' which became known as the Littlewood-Parley conjecture. The truth of this conjecture for certain subclasses of $\mathcal{S}$ enshrouded its falsity in general until as early as in 1933 (about a year after the conjecture), when it was settled in negation by what came to be known as the Fekete-Szeg$\ddot{o}$ problem.\vskip 2mm

{\bf Fekete-Szeg$\ddot{o}$ problem}. The origin of this problem is the disproof of the conjecture of Littlewood and Parley with regard to the bound on the coefficient of odd univalent functions as has preceeded. For each $f\in\mathcal{S}$, Fekete and Szeg$\ddot{o}$ obtained the sharp bound:
$$|a_3-\alpha a_2^2|\leq 1+2e^{-2\alpha/(1-\alpha)},\;\;0\leq\alpha\leq 1.$$
This results gives $|c_5|<1/2+e^{-2/3}=1.013\cdots$ because $c_5=(a_3-a_2^2/4)/2$.\vskip 2mm

Thus the Fekete-Szeg$\ddot{o}$ problem has continued to recieve attention until even in the many subclasses of $\mathcal{S}$. The functional $|a_3-\alpha a_2^2|$ is well known as the Fekete-Szeg$\ddot{o}$ functional. Many other functionals have risen after it, each finding application in certain problems of the geometric functions. For $\alpha=1$, it is important to mention a more general problem of this type, which is the Hankel determinant problem.\vskip 2mm

{\bf Hankel determinant problem}. Let $n\geq 0$ and $q\geq 1$, the $q$-th Hankel determinant of the coefficients of $f\in\mathcal{S}$ is defined as:
\[
H_q(n)=
\begin{vmatrix}
a_n & a_{n+1} & \cdots & a_{n+q-1}\\
a_{n+1} & \cdots & \cdots & \vdots\\
\vdots & \vdots & \vdots & \vdots\\
a_{n+q-1} & \cdots & \cdots & a_{n+2(q-1)}
\end{vmatrix}.
\]
The determinant has been investigated by several authors with the subject of inquiry ranging from rate of growth of $H_q(n)$ as $n\rightarrow\infty$ to the determination of precise bounds on $H_q(n)$ for specific $q$ and $n$ for some favored classes of functions. It is interesting to note that $|H_2(1)|\equiv|a_3-a_2^2|$, the Fekete-Szeg$\ddot{o}$ functional for $\alpha=1$.\vskip 2mm

{\bf Other coefficient related problems}. These include the determination of successive coefficient relationship and the region of variability of coefficients.\vskip 2mm

{\bf Growth, Distortion and Covering}. The idea of growth of analytic function $f$ refers to the size of the image domain, that is $|f(z)|$. The term, distortion, arises from the geometric interpretation of $|f'(z)|$ as the infinitesimal magnification factor of the arclength under the mapping $f$, or from the Jacobian $|f'(z)|^2$ as the infinitesimal magnification factor of the area of the image domain. The concept of covering by a function $f$ refers to the portion of the image domain covered by it. For the large family of univalent functions, it is known that the range of every member function covers the disk $|\xi|<1/4$. This assertion is due to Koebe, 1907, and has thus been known as the Koebe One-Quarter Theorem. It is a consequence of the Bieberbach Theorem on the second coefficient of functions in $S$ and their omitted-value transformation.\vskip 2mm

\begin{proof}
If $f\in\mathcal{S}$ omits $\xi\in\mathbb{C}$, then
$$g(z)=\frac{\xi f(z)}{\xi-f(z)}=z+\left(a_2+\frac{1}{\xi}\right)z^2\cdots$$
is analytic and univalent in $E$. So by Bieberbach theorem
$$\left|a_2+\frac{1}{\xi}\right|\leq 2$$
combined with the fact that $|a_2|\leq 2$, the covering $|\xi|<1/4$ follows.
\end{proof}
\vskip 2mm

{\bf Partial sums}. The inquistion regarding partial sums
$$s_n(z)=z+a_2z^2+\cdots+a_nz^n$$
of the series development of $f$ is about the extent to which known geometric properties of $f$ are carried on to its partial sums. Another result of Yamaguchi[] is suitable to mention here:\vskip 2mm

\begin{quote}[Yamaguchi \cite{KY}]
If $f$ satisfies Re $f(z)/z>0$ in $E$, then the $k$th partial sums $s_k(z)=z+a_2z^2+\cdots+a_kz^k$ is univalent in the subdisk $|z|<\frac{1}{4}.$
\end{quote}
\vskip 2mm

{\bf Linear sums or combinations}. It is also of interest to find out: if $\phi$ and $\varphi$ are some geometric quantities about $f$, then under what conditions is the linear sum $(1-t)\phi+t\varphi$ preserving some known geometric properties based on $\phi$ and $\varphi$?

\section{Some subclasses of $\mathcal{S}$}
Sequel to what has preceeded of some of the subclasses of the class of univalent functions, we mention that the fundamental basis or justification for discussing new subclasses lies in the fact through them certain classes of functions may be associated with some special properties, not commonly associable with certain other classes. Thus the many subjects of inquiry are being reinvestigated in several class of functions to sharpen, smoothen or better many known results particularly in the direction of a new subclass. Some of the well known subclasses of $S$ (with the associated geometric quantities in brackets) are:\vskip 2mm

{\bf Functions of bounded turning ($f'$).} These are functions whose derivatives have positive real parts, that Re $f'(z)>0$. They are entirely univalent functions as has preceeded. Many results concerning this can be found in the literatures.\vskip 2mm

{\bf Starlike functions ($zf'/f$).} They are functions for which the real part of the quantity $zf'/f$ is positive. They are entirely univalent functions. They are also convex. They are close-to-convex as well. Results on this class of functions are scattered in many literatures.\vskip 2mm 

{\bf Convex functions ($1+zf''/f'$).} They are functions for which the real part of the quantity $1+zf''/f'$ is positive. They are entirely univalent functions. They are also close-to-convex. Results on this class of functions are scattered in many literatures.\vskip 2mm

{\bf Quasi-convex ($(zf')'/g'$, $g$ is convex).} They are functions for which the real part of the quantity $(zf')'/g'$, $g$ is convex, is positive. They are entirely univalent functions. They are also close-to-convex. Results in this direction are also scattered in many literatures.\vskip 2mm

{\bf Close-to-convex ($f'/g'$, $g$ is convex).} They are functions for which the real part of the quantity $f'/g'$, $g$ is convex, is positive. They are entirely univalent functions. Results in this direction are also scattered in many literatures. A bounded turning function is a special close-to-convex function with $g(z)=z$.\vskip 2mm

{\bf Bazilevic functions.} They consist of functions defined by the integral
$$f(z)=\left\{\frac{\alpha}{1+\xi^2}\int_0^z[h(t)-i\xi]t^{-\left(1+\frac{i\alpha\xi}{1+\xi^2}\right)}g(t)^{\left(\frac{\alpha}{1+\xi^2}\right)}dt\right\}^{\frac{1+i\xi}{\alpha}}$$
where $h$ is an analytic function which has positive real part in $E$ and normalized by $h(0) = 1$ and $g$ is starlike  in $E$.  The numbers $\alpha>0$ and $\xi$ are real and all powers meaning principal determinations only. They are entirely univalent in the unit disk. They contain many other class of function as special cases.\vskip 2mm

{\bf Inclusions.} Two well known inclusion relations between these classes are given as:\vskip 2mm

{\em convexity $\Longrightarrow$ quasi-convexity $\Longrightarrow$ close-to-convexity $\Longrightarrow$ univalence}.\vskip 2mm 

{\em convexity $\Longrightarrow$ starlikeness $\Longrightarrow$ close-to-convexity $\Longrightarrow$ univalence}.\vskip 2mm 

{\bf Other subclasses and generalizations.} There are many other subclasses of the above classes of functions which have appeared in prints. Many generalizations have also appeared via derivative as well as integral operators. These operators include the well known Salagean derivative, Ruscheweyh derivative, Noor integral operator and some further generalizations of them.

\section{Caratheodory, related functions and generalizations}
A cursory look at the series development (1) for $f$ and the various geometric quantities $zf'/f$,  $1+zf''/f'$, $f/g$, $f'/g'$, and many more, (which possess the property of positivity of real parts) suggests clearly the existence of a series form:
$$h(z)=1+c_1z+c_2z^2+\cdots.\eqno(2)$$
The form (2) satisfies $h(0)=1$ and Re $h(z)>0$ (positive real parts). The present author is not aware the discovery of which predates which of the two functions, $f$ (normalized by $f(0)=0$ and $f'(0)=1$) and $h$ (normalized by $h(0)=1$). However, it is not out of place to insinuate that the prediscovery $f$ over $h$. For otherwise, the discovery of $f$ certainly would have spurred inquisition into $h$. The study of $h$ provides much insight into the natures of any $f$ having goemetries described above. The function $h$ is called the Caratheodory function (named after Caratheodory who not only noticed the obvious, but expended much energy in its characterizations). The function $h$ may be described equivalently as a function {\em subordinate} to the M$\ddot{o}$ebius function,
$$L_0(z)=\frac{1+z}{1-z}.$$

The M$\ddot{o}$ebius function play a central role in the family of functions of the like of $h$. It assumes the extremum in the most extremal problem for such functions.\vskip 2mm

By {\em subordination}, it is meant that there exist a function of unit bound, $\vartheta(z)$ ($|\vartheta(z)|<1$, normalized by $\vartheta(0)=0$) such that $h(z)=L_0(\vartheta(z))$. Thus this gives another representation for $h$ among others.  Precisely, in terms $\vartheta$, $h$ has the form:
$$h(z)=\frac{1+\vartheta(z)}{1-\vartheta(z)},\;\; z\in E.$$

The unit bound functions are known as Schwarz functions. Two basic results are noteworthy about them. These are:

\begin{quote}[Schwarz (See \cite{CC})]
{\em If $\vartheta(z)$ is a function of unit bound in $E$, then for each $0<r<1$, $|\vartheta(0)|<1$ and $|\vartheta(re^{i\theta})|\leq r$ unless $\vartheta(z)=e^{i\sigma}z$ for some real number $\sigma$.} 
\end{quote}

The above result is commonly refered to as the Schwarz's Lemma. It has the implication that if $\vartheta(z)$ is a function of unit bound in $E$, so also is $u(z)=\vartheta(z)/z$, that is $|u(z)|<1$, but not necessarily normalized by $|u(0)|=0$.

\begin{quote}[Caratheodory \cite{CC}]
{\em If $\vartheta(z)$ is a function of unit bound (not necessarily normalized) in $E$, then
$$|\vartheta'(z)|\leq\frac{1-|\vartheta(z)|^2}{1-|z|^2}$$
with strict inequality holding unless $\vartheta(z)=e^{i\sigma}z$ for some real number $\sigma$.}
\end{quote}

Studies have also revealed that any $h$ can as well have what is known as the Herglotz representation, which is the integral form:
$$h(z)=\int_0^{2\pi}\frac{e^{it}+z}{e^{it}-z}d\mu(t),$$
where $d\mu(t)\geq 0$ and $\int d\mu(t)=\mu(2\pi)-\mu(0)=1$. \vskip 2mm

The various represntations of $h$ have important applications as may be discovered through further studies.\vskip 2mm

The Caratheodory functions are also preserved under a number of transformations: suppose $g,h$ are Caratheodory, then so is $p$ defined as (i) $p(z)=g(e^{it}z)$, $t$ real; (ii) $p(z)=g(tz)$, $t\in[-1,1]$; (iii) $p(z)=g[(z+t)/(1+\bar{t}z)]/g(t)$, $|t|<1$; (iv) $p(z)=(g(z)+it)/(1+itg(z))$, $t$ real; (v) $p(z)=[g(z)]^t$, $t\in[-1,1]$ and (vi) $p(z)=[g(z)]^t[h(z)]^\tau$, $t,\tau,t+\tau\in[0,1]$.\vskip 2mm

\begin{proof}
By simple computation it is easy to see that in all cases, $p(0)=1$. Thus it only remains to show that the real parts of the transformations are positive. For (i) - (iii), this follows from the fact that each of the points $e^{it}z$, $tz$ and $(z+t)/(1+\bar{t}z)$, (with associated conditions on $t$) are transformations of points in $|z|<1$ to points in there$^{(6)}$\footnote{ $^{(6)}$To show that $|(z+t)/(1+\bar{t}z)|<1$, assume the converse. That is $|z+t|\geq|1+\bar{t}z|$. Then squaring both sides we obtain $|z|^2+|t|^2\geq 1+|t|^2|z|^2$, wherefrom we obtain a contradiction that $|z|\geq 1$. This proves the point.}. In fact (iv) is a linear transformation of the right half plane onto itself$^{(7)}$\footnote{ $^{(7)}$The fact that $p(z)=(g(z)+it)/(1+itg(z))$ is a linear transformation of the right half plane onto itself can be deduced from the fact that:
$Re\left\{\frac{g+it}{1+itg}\right\}=Re\left\{\frac{(g+it)(1-it\bar{g})}{(1+itg)(1-it\bar{g})}\right\}=\frac{Re(g+t^2\bar{g})}{|1+itg|^2}>0$ since the real parts of $g$ and $\bar{g}$ is greater that zero.} while (v) and (vi) follow from the fact that $Re\;z^t\geq (Re\;z)^t$ when $t\in[0,1]$ and $Re\;z>0^{(8)}$\footnote{ $^{(8)}$The fact that $Re\;z^t\geq (Re\;z)^t$ is due, by elementary calculus, to the fact that $y=\cos t\theta/(\cos{\theta})^t$ attains its maximum value at $t_0\in[0,1]$ ($t_0$ is given by $t_0=\frac{\arctan\left(\frac{-\log\cos\theta}{\theta}\right)}{\theta}$, for all $\theta\in(-\frac{\pi}{2},\frac{\pi}{2})$, $\theta\neq 0$) and $y(t)$ is decreasing on $t\in[t_0,1]$; and $y(t)$ is increasing on $t\in[0,t_0]$. In particular, $y(t)=\cos t\theta/(\cos{\theta})^t\geq y(0)=y(1)=1$.}. Then for each $t\in[-1,0]$ with respect to (v), the function $p$ takes the reciprocal of its values for $t\in[0,1]$. This concludes the proof.
\end{proof}

We now mention two basic coefficient inequalities for $h$, the first based on its Herglotz representation while the other depends on its representation by functions of unit bound $\vartheta(z)$.

\begin{quote}[Caratheodory (See \cite{PD})]
{\em If $h(z)=1+c_1z+c_2z^2+\cdots$ is a Caratheodory function, then $|c_k|\leq 2$, $k=1,2,\cdots$. The M$\ddot{o}$ebius function takes the equality.} 
\end{quote}

\begin{proof}
If we expand the Herglotz representation of $h$ in series form$^{(9)}$\footnote{ $^{(9)}$The Herglotz representation $h(z)=\int_0^{2\pi}\frac{e^{it}+z}{e^{it}-z}d\mu(t)$
can be written as $$\aligned h(z)&=\int_0^{2\pi}\frac{1+ze^{-it}}{1-ze^{-it}}d\mu(t)\\ &=\int_0^{2\pi}\left(1+2ze^{-it}+2z^2e^{-2it}+2z^3e^{-3it}+\cdots\right)d\mu(t)\\ &=1+\sum_{k=1}^\infty\left(2\int_0^{2\pi}e^{-ikt}d\mu(t)\right)z^k.\endaligned$$ So, when compared with $h(z)=1+\sum_{k=1}^\infty c_kz^k$ gives $c_k=2\int_0^{2\pi}e^{-ikt}d\mu(t)$.} and compare coefficients of $z^k$, we find that;
$$c_k=2\int_0^{2\pi}e^{-ikt}d\mu(t).$$
So, we have
$$\aligned
|c_k|
&\leq 2\int_0^{2\pi}|e^{-ikt}|d\mu(t),\;\text{since}\;d\mu(t)\;\text{is nonnegative}\\
&=2\int_0^{2\pi}d\mu(t)=2,\;\text{since}\;\int_0^{2\pi}d\mu(t)=1.\endaligned$$
\end{proof}

\begin{quote}[Pommerenke \cite{CP}]
{\em If $h(z)=1+c_1z+c_2z^2+\cdots$ is a Caratheodory function, then
$$\left|c_2-\frac{c_1^2}{2}\right|\leq 2-\frac{|c_1|^2}{2}.$$}
Equality holds for the function:
$$h(z)=\frac{1+\frac{1}{2}(c_1+\varepsilon\overline{c_1})z+\varepsilon z^2}{1-\frac{1}{2}(c_1-\varepsilon\overline{c_1})z-\varepsilon z^2},\;|\varepsilon|=1.$$ 
\end{quote}

\begin{proof}
Suppose $\vartheta(z)$ is a function of unit bound in $E$, normalized by $\vartheta(0)=0$. Then by Schwarz's Lemma there exists an analytic function $u(z)$ also of unit bound such that $\vartheta(z)=zu(z)$. Then
$$\aligned
u(z)=\vartheta(z)/z
&=\frac{1}{z}\frac{h(z)-1}{h(z)+1}\\
&=\frac{1}{2}c_1+\left(\frac{1}{2}c_2-\frac{1}{4}c_1^2\right)z+\cdots\endaligned$$ 
satisfies $|u(z)|\leq 1$ in $E$ so that
$$|u'(z)|\leq\frac{1-|u(z)|^2}{1-|z|^2}.$$
Thus
$$\aligned
|u'(0)|
&=\frac{1}{2}c_2-\frac{1}{4}c_1^2\\
&\leq 1-|u(0)|^2=1-\frac{1}{4}|c_1|^2.\endaligned$$
\end{proof}

The above two basic inequalities have great implications in the field, especially with regard to coefficient problems.\vskip 2mm

Further advances have led to various generalizations of $h$. Janowski \cite{JW} redefined $h$ in terms of $\vartheta$, saying given fixed real numbers $a,b$ such that $a\in(-1,1]$ and $b\in[-1,a)$ (that is $-1\leq b<a\leq 1$), then $h$ is defined as:
$$h(z)=\frac{1+a\vartheta(z)}{1+b\vartheta(z)},$$
where the Caratheodory function corresponds to the extremes $b=-1$, $a=1$. For various choice values of $a,b$, the function $h$ also maps the unit disk to some portions of the right half plane.\vskip 2mm

Perhaps, if any, the most significant of our contribution to this important field of study is the development of iterations for the very important families of functions: the Caratheodory and Janowski functions (See \cite{KO} and the previous works cited therein). These are:
$$p_n(z)=\frac{\alpha}{z^\alpha}\int_0^zt^{\alpha-1}p_{n-1}(t)dt,\;\;n\geq 1,$$
with $p_0(z)=p(z)$.
$$p_{\sigma,n}(z)=\frac{\sigma-(n-1)}{z^{\sigma-(n-1)}}\int_0^zt^{\sigma-n}p_{\sigma,n-1}(t)dt,\;\;\;
n\geq 1$$
with $p_{\sigma,0}(z)=p(z)$.
These transformations preserve many geometric structures of the family of functions with positive real part normalized by $h(0)=1$; particularly the positivity of the real parts, compactness, convexity and subordination. Another fascinating aspect of these transformations is that with them, investigations of the various classes of functions associated with them have become easy, short and elegant. They have been very helpful in dealing easily with certain problems of the theory of analytic and univalent function in the most intriguing simplicity.\vskip 2mm

Many techniques have been developed in the field. However, the most fundamental and beginner-friendly is one based on the close association that exists between the Caratheodory functions (together with its further developments) and many classes of functions. Many fundamental results have been established as regards this class of functions. Thus investigating various problems of geometric functions via an underlying $h$ has been well accepted among researchers in this field as a princpal technique. In the next section are given a few examples and insight into the technique of constructing the extremal functions.

\section{Illustrating examples}
The examples in this section are simple. Our objective is to demonstrate, in repeated and beginner-friendly manner, how results can be obtained in certain classes of functions using an underlying Caratheodory function, $h$. All the results here are best possible. The construction of extremal functions are greatly simplified.\vskip 2mm

{\bf Theorem A.} {\em If $f\in\mathcal{S}$ satisfy Re $f(z)/z>0$, then its coefficients satisfy the inequality: $|a_k|\leq 2$. Equality is attained by $f(z)=z(1+z)/(1-z)$.}

\begin{proof}
Since Re $f(z)/z>0$, then $f(z)/z$ is a function with positive real parts. Hence, $f(z)/z=h(z)$ for some $h(z)$ with positive real parts. Comparing coefficients of the series expansion of $f(z)/z$ and $h$, we have $a_k=c_{k-1}$, $k=2,3,\cdots$ so that $|a_k|\leq 2$ since $|c_k|\leq 2$, $k=1,2,\cdots$.\vskip 2mm

The construction of the extremal function is by setting the geometric quantity $f(z)/z$ equal to the extremal function for $h$, which is $L_0(z)=(1+z)/(1-z)$. This simply gives $f(z)=z(1+z)/(1-z)$.
\end{proof}

{\bf Theorem B.} {\em The coefficients of functions of bounded turning $($Re $f'(z)>0)$ satisfy the inequality: $|a_k|\leq 2/k$. Equality is attained by $f(z)=-2\ln(1-z)-z$.}

\begin{proof}
Since Re $f'(z)>0$, then $f'$ is a function with positive real parts. Hence, $f'(z)=h(z)$ for some $h(z)$ with positive real parts. Comparing coefficients of the series expansion of $f'$ and $h$, we have $a_k=c_{k-1}/k$, $k=2,3,\cdots$ so that $|a_k|\leq 2/k$ since $|c_k|\leq 2$, $k=1,2,\cdots$.\vskip 2mm

The construction of the extremal function is by setting the geometric quantity $f'(z)$ equal to the extremal function for $h$, which is $L_0(z)=(1+z)/(1-z)$. Thus,  we have
$$f'(z)=\frac{1+z}{1-z}.$$
Integrating both sides, we have
$$f(z)=\int_0^z\frac{1+t}{1-t}dt=-\ln(1-z)-z.$$
\end{proof}

{\bf Theorem C.} {\em If $f$ is starlike $($Re $zf'(z)/f(z)>0)$, then $|a_k|\leq k$. Equality is attained by the Koebe function $k(z)=z/(1-z)^2$.}

\begin{proof}
Given that $f(z)=z+a_2z^2+\cdots$. Since Re $zf'(z)/f(z)>0$, then $zf'(z)/f(z)$ is a function with positive real parts. Hence, $zf'(z)/f(z)=h(z)$ for some $h(z)=1+c_1z+\cdots$ with positive real parts. Equating $zf'(z)/f(z)$ and $h$ we have
$zf'(z)=h(z)f(z)$. Expanding both sides in series, we have
$$\aligned
zf'(z)
&=z+2a_2z^2+3a_3z^3+\cdots\\
&=f(z)h(z)\\
&=z+(a_2+c_1)z^2+(a_3+a_2c_1+c_2)z^3+(a_4+a_3c_1+a_2c_2+c_3)z^4+\cdots\endaligned$$
so that
$$ka_k=a_k+\sum_{j=1}^{k-1}a_jc_{k-j},\;a_1=1$$
and thus
$$(k-1)a_k=\sum_{j=1}^{k-1}a_jc_{k-j},\;a_1=1.$$
We now proceed by induction. For $k=2$, we have $a_2=a_1c_1$ with $a_1=1$ so that $|a_2|=|c_1|\leq 2$ as required. Next we suppose the inequality is true for $k=n$, then for $k=n+1$ we have
$$na_{n+1}=\sum_{j=1}^{n}a_jc_{n+1-j}$$
so that 
$$n|a_{n+1}|\leq\sum_{j=1}^{n}|a_j||c_{n+1-j}|\leq 2\sum_{j=1}^{n}j=n(n+1).$$
Thus we have $|a_{n+1}|\leq n+1$ and the inequality follows by induction.\vskip 2mm

As for the extremal function, the construction is by setting the geometric quantity $zf'(z)/f(z)$ equal to the extremal function for $h$, which is $L_0(z)=(1+z)/(1-z)$. Thus, we have
$$\frac{zf'(z)}{f(z)}=\frac{1+z}{1-z}$$
so that
$$\frac{f'(z)}{f(z)}=\frac{1+z}{z(1-z)}=\frac{1}{z}+\frac{2}{1-z}.$$
Now integrating both sides, we have $\ln f(z)=\ln z-2\ln(1-z)$ which gives $f(z)=z/(1-z)^2$, which is the Koebe function.
\end{proof}

{\bf Theorem D.} {\em If $f$ is convex $($Re $[1+zf''(z)/f'(z)]>0)$, then $|a_k|\leq 1$. Equality is attained by the Koebe function $f(z)=1/(1-z)$.}

\begin{proof}
Observe that we can write the geometric quantity $1+zf''(z)/f'(z)$ as $(zf'(z))'/f'(z)$ so that the convexity condition Re $[1+zf''(z)/f'(z)]>0$ now becomes Re $[z(zf'(z))']/[zf'(z)]>0$. That is $zf'(z)$ is starlike (in fact $f$ is convex if and only if $zf'(z)$ is starlike)$^{(10)}$\footnote{ $^{(10)}$The statement ``$f$ is convex if and only if $zf'(z)$ is starlike'' is due to Alexander and is known as Alexander theorem}. So by the result for starlike function, the coefficients of $zf'(z)$ satisfy $|a_k|\leq k$. Hence we have $k|a_k|\leq k$ which gives $|a_k|\leq 1$ as required.\vskip 2mm

Now to the construction of the extremal function, set the geometric quantity $1+zf''(z)/f'(z)$ equal to the extremal function for $h$, which is $L_0(z)=(1+z)/(1-z)$. Thus, we have
$$\frac{z(zf'(z))'}{zf'(z)}=\frac{1+z}{1-z}$$
so that as in the previous proof we have $zf'(z)=z/(1-z)^2$. Furthermore we have $f'(z)=1/(1-z)^2$, which on integration gives $f(z)=1/(1-z)$.
\end{proof}

{\bf Theorem E.} {\em If $f$ is close-to-convex $($Re $f'(z)/g'(z)>0$, $g$ is convex$)$, then $|a_k|\leq k$. Equality is attained by the Koebe function $k(z)=z/(1-z)^2$.}

\begin{proof}
Given that $f(z)=z+a_2z^2+\cdots$. Since Re $f'(z)/g'(z)>0$, then $f'(z)/g'(z)$ is a function with positive real parts. Hence, $f'(z)/g'(z)=h(z)$ for some $h(z)=1+c_1z+\cdots$ with positive real parts. Equating $f'(z)/g'(z)$ and $h$ we have
$f'(z)=h(z)g'(z)$. Let $g(z)=z+b_2z^2+\cdots$. Expanding both sides in series, we have
$$\aligned
f'(z)
&=1+2a_2z+3a_3z^2+\cdots\\
&=h(z)g'(z)\\
&=1+(2b_2+c_1)z+(3b_3+2b_2c_1+c_2)z^2+(4b_4+3b_3c_1+2b_2c_2+c_3)z^3+\cdots\endaligned$$
so that
$$ka_k=kb_k+\sum_{j=1}^{k-1}jb_jc_{k-j},\;b_1=1$$
and thus
$$k|a_k|\leq k|b_k|+\sum_{j=1}^{k-1}j|b_j||c_{k-j}|,\;a_1=1.$$
Since $|b_k|\leq 1$, $k=2,3,\cdots$ for convex functions and $|c_k|\leq 2$, $k=1,2,\cdots$ for $h$, it follows that $k|a_k|\leq k+2\sum_{j=1}^{k-1}j=k+k(k-1)$ so that the desired inequality follows.\vskip 2mm

As for the extremal function, the construction is by choosing $g=1/(1-z)$ in the geometric quantity $f'(z)/g'(z)$ and setting this equal to the extremal function for $h$, which is $L_0(z)=(1+z)/(1-z)$. Thus, we have
$$\frac{f'(z)}{g'(z)}=\frac{f'(z)}{1/(1-z)^2}=\frac{1+z}{1-z}$$
so that
$$f'(z)=\frac{1+z}{(1-z)^3}.$$
Integrating both sides, we have $f(z)=z/(1-z)^2$, which is the Koebe function.
\end{proof}

\medskip

{\it Acknowledgements.} The author is indebted to his undergraduate project students: Balqees Alege, Charles Nwawulu, Aminah Olasupo and Olawale Osewa all at the University of Ilorin, Ilorin, Nigeria, whose interest in the subject informed the compilation of this work.

\bigskip

\bibliographystyle{amsplain}

\begin{thebibliography}{99}

\bibitem {OP}
\textsc{Ahuja, O. P.}, {\it Planar harmonic univalent and related mappings}, J. Inequal. Pure and Applied Math., {\bf6} (4) Art. 122 (2005), 1--18.

\bibitem {KO}
\textsc{Babalola, K. O.}, {\it Convex null sequence technique for analytic and univalent mappings of the unit disk}, Tamkang J. Math., {\bf40} (2) (Summer 2009), 201--209.

\bibitem {IE}
\textsc{Bazilevic, I. E.}, {\it On a case of integrability by quadratures of the equation of Loewner-Kufarev} Mat. Sb. {\bf37} (79), (1955), 471--476. (Russian).

\bibitem {CC}
\textsc{Caratheodory, C.}, {\it Theory of functions of a complex variable}, Vols. I \& II, Chelsea Pub. Co. New York Inc. 1960.

\bibitem {PL}
\textsc{Duren, P. L.}, {\it Coefficients of univalent functions}, Bull. Amer. Math. Soc. {\bf83} (1977), 891--911.

\bibitem {PD}
\textsc{Duren, P. L.}, {\it Univalent functions}, Springer Verlag. New York Inc. 1983.

\bibitem {AW}
\textsc{Goodman, A. W.}, {\it Univalent Functions}, Vols. I \& II, Marina Pub. Co., Tampa, Florida, 1993.

\bibitem {JW}
\textsc{Janowski, W.}, {\it Some extremal problems for certain families of analytic functions I}, Ann. Polon. Math. {\bf28}, (1973), 297--326.

\bibitem {RJ}
\textsc{Libera, R. J.}, {\it Some classes of regular univalent functions}, Proc. Amer. Math. Soc. {\bf16} (4) (1965), 755--758.

\bibitem {TH}
\textsc{Macgegor, T. H.}, {\it Functions whose derivative has a positive real part}, Trans. Amer. Math. Soc. {\bf104} (1962), 532--537.

\bibitem {TM}
\textsc{Macgegor, T. H.}, {\it Geometric problems in complex analysis}, The Amer. Math. Monthly {\bf79} (5) (1972), 447--468.

\bibitem {CP}
\textsc{Pommerenke, Ch.}, {\it Univalent Functions}, Vandenhoeck and Ruprecht, G$\ddot{o}$ttingen, 1975.

\bibitem {KY}
\textsc{Yamaguchi, K.}, {\it On functions satisfying $Re\{f(z)/z\}>0$}, Proc. Amer. Math. Soc. {\bf17} (1966), 588--591.

\end{thebibliography}

\end{document}